\documentclass[12pt]{article}
\usepackage[utf8]{inputenc}
\usepackage{setspace}
\usepackage{graphicx}
\usepackage{array}
\usepackage[dvipsnames]{xcolor}
\usepackage{amsmath}
\usepackage{fancyvrb}
\usepackage{mathtools}
\usepackage{amssymb}
\usepackage{comment}
\usepackage{float}
\usepackage{booktabs}
\usepackage[bottom]{footmisc}
\usepackage[shortlabels]{enumitem}
\usepackage{tikz-cd}
\usepackage{hyperref}
\usepackage{cleveref}
\usepackage{amsthm}
\usepackage{xfrac}
\usepackage[english]{babel}
\usepackage{url}
\usepackage{yfonts}

\usepackage[margin=0.8in]{geometry}

\newcommand{\C}{\ensuremath{\mathbb{C}}}
\newcommand{\F}{\ensuremath{\mathbb{F}}}
\newcommand{\Z}{\ensuremath{\mathbb{Z}}}

\hypersetup{
    colorlinks=true, 
    linktoc=all,     
    linkcolor=black,  
    allcolors=WildStrawberry
}

\newcommand{\thistheoremname}{}
\newtheorem*{genericthm*}{\thistheoremname}
\newenvironment{namedthm*}[1]
  {\renewcommand{\thistheoremname}{#1}%
   \begin{genericthm*}}
  {\end{genericthm*}}

\usepackage{mathabx}

\usepackage{mathtools}
\usepackage{mathdots}
\usepackage[document]{ragged2e}

\theoremstyle{definition}

\theoremstyle{definition}

\theoremstyle{definition}
\newtheorem{theorem}{Theorem}

\theoremstyle{definition}

\theoremstyle{definition}

\theoremstyle{definition}

\theoremstyle{definition}

\theoremstyle{definition}

\theoremstyle{definition}

\theoremstyle{definition}
\newtheorem{proposition}{Proposition}

\title{\textbf{On representation zeta function of special linear groups over finite principal ideal local rings}}
\author{Uri Ronen\\
\href{mailto:ur@campus.technion.ac.il}{\color{black}ur@campus.technion.ac.il}}
\date{\vspace{-5ex}}

\begin{document}

\maketitle

\begin{abstract} 
    We show that the group algebras $\C[\text{SL}_3(\sfrac{\F_3[t]}{(t^3)})]$ and $\C[\text{SL}_3(\sfrac{\Z}{27})]$ are not isomorphic, as well as $\C[\text{SL}_4(\sfrac{\F_2[t]}{(t^3)})]\ncong\C[\text{SL}_4(\sfrac{\Z}{8})]$, by computing the number of conjugacy classes in those groups using MAGMA's calculator. Similarly, we reproduce special cases of a recent result by Hassain and Singla, showing that $\C[\text{SL}_2(\sfrac{\F_2[t]}{(t^k)})]\ncong\C[\text{SL}_2(\sfrac{\Z}{2^k})]$ for $3\leq k\leq 8$.
\end{abstract}
\begin{center}
\section*{Introduction}
\end{center}
For a finite group $G$, let $\text{Irr}(G)$ be the set of complex irreducible representations of $G$, up to equivalence. Define the representation zeta function of $G$ to be:
\begin{equation*}  \zeta_G(s):=\sum_{\rho\in\text{Irr}(G)} \frac{1}{(\dim\rho)^s}
\end{equation*} 
Note that if $G_1,G_2$ are finite groups then $\zeta_{G_1}=\zeta_{G_2}$ if and only if their groups algebras over $\C$ are isomorphic. Let $p$ be a prime and let $\Z_p$ be the ring of $p$-adic integers, and $\F_p[[t]]$ the ring of power series over $\F_p$. When studying (continuous) representations of profinite groups such as $\text{GL}_n(\F_p[[t]])$ and $\text{GL}_n(\Z_p)$, it is natural to study the representations of their finite quotients, as continuous representations of profinite groups factor through some finite quotient. In this context, we have the following (special case of a) conjecture by Onn (\cite{ZetaFunctionOfQuotientsOfGL2}, conjecture 1.3):
\begin{equation*}
    \zeta_{\text{GL}_n(\sfrac{\F_p[t]}{(t^k)})}=\zeta_{\text{GL}_n(\sfrac{\Z}{p^k})}\text{ for all }n,k\geq1 \text{ and }p \text{ prime}.
\end{equation*}
As of the time of writing this paper, Onn's conjecture remains open. Nevertheless, extensive efforts have been devoted over the years to establish it's correctness: in \cite{ZetaFunctionOfQuotientsOfGL2} the conjecture was demonstrated to hold for $n=2$ and all $k\geq1$ and $p$. Recently, \cite{ZetaFunctionOfSL3WithoutRestrictionOnPrime} extended this result to $\text{SL}_3$ and $\text{GL}_3$, under the condition that $p>3$. Subsequent work in \cite{ZetaFunctionOfLengthTwoRingsOfGL} proved Onn’s conjecture for $k=2$ and for all $n\geq1$ and $p$, and shortly after \cite{ZetaFunctionOfLengthTwoRingsOfClassicalGroups} further extended this result to other classical groups, such as $\text{SL}_n,\text{O}_n,\text{Sp}_n,\text{U}_n$. A recent contribution in \cite{LogicResultByItamar} established that for every $n,k\geq1$ there exists $N(n,k)>0$ such that for any prime $p>N(n,k)$ Onn's conjecture holds with $n,k,p$, and the same is true for $\text{SL}$. \\
Conversely, in (\cite{CounterexampleForSL2(EvenDegrees)},\cite{CounterexampleForSL2(OddDegrees)}) it was proven that the conjecture is false for $\text{SL}_2$, $p=2$, and $k\geq4$. Notably, this result can be extended to   $k=3$ using similar methods, but not to $k=2$, as demonstrated in \cite{CounterexampleForSL2(EvenDegrees)}  ,\S11.\\
With the assistance of MAGMA's online computational algebra calculator, we found explicit generators for some of the groups mentioned above and subsequently determined their number of conjugacy classes: 
\begin{center}
\begin{tabular}{|c| c!{\vrule width 2pt}c |c|}
    \hline
    $G_1$ &\multicolumn{2}{|c|}{\#\text{Conj.Classes}} &$G_2$ \\
    \hline
    $\text{SL}_4(\sfrac{\F_2[t]}{(t^3)})$ &1824 &1896&$\text{SL}_4(\sfrac{\Z}{2^3})$\\
    \hline
    $\text{SL}_3(\sfrac{\F_3[t]}{(t^3)})$ &1242&1218&$\text{SL}_3(\sfrac{\Z}{3^3})$\\
    \hline
    $\text{SL}_2(\sfrac{\F_2[t]}{(t^8)})$&1336&1456&$\text{SL}_2(\sfrac{\Z}{2^8})$\\
    \hline
    $\text{SL}_2(\sfrac{\F_2[t]}{(t^7)})$&624&720&$\text{SL}_2(\sfrac{\Z}{2^7})$\\
    \hline
    $\text{SL}_2(\sfrac{\F_2[t]}{(t^6)})$&292&352&$\text{SL}_2(\sfrac{\Z}{2^6})$\\
    \hline
    $\text{SL}_2(\sfrac{\F_2[t]}{(t^5)})$&132&168&$\text{SL}_2(\sfrac{\Z}{2^5})$\\
    \hline
    $\text{SL}_2(\sfrac{\F_2[t]}{(t^4)})$&58&76&$\text{SL}_2(\sfrac{\Z}{2^4})$\\
    \hline
    $\text{SL}_2(\sfrac{\F_2[t]}{(t^3)})$&24&30&$\text{SL}_2(\sfrac{\Z}{2^3})$\\
    \hline
\end{tabular}
\end{center}
Note that if $\zeta_{G_1}=\zeta_{G_2}$,  then in particular $\#\text{Conj.Classes}(G_1)=\#\text{Irr}(G_1)=\zeta_{G_1}(0)=\zeta_{G_2}(0)=\#\text{Irr}(G_2)=\#\text{Conj.Classes}(G_2)$, so it follows from the above:
\begin{theorem}
    The following pairs of groups have nonisomorphic complex group algebra, or equivalently, different representation zeta function:
    \begin{itemize}
        \item $\text{SL}_4(\sfrac{\F_2[t]}{(t^3)})$ and $\text{SL}_4(\sfrac{\Z}{2^3})$
        \item $\text{SL}_3(\sfrac{\F_3[t]}{(t^3)})$ and $\text{SL}_3(\sfrac{\Z}{3^3})$
    \end{itemize}
\end{theorem}
Additionally, this data reproduces special cases of the main result presented in (\cite{CounterexampleForSL2(EvenDegrees)},\cite{CounterexampleForSL2(OddDegrees)}):
\begin{proposition} $\C[\text{SL}_2(\sfrac{\F_2[t]}{(t^k)})]\ncong\C[\text{SL}_2(\sfrac{\Z}{2^k})]$ for $3\leq k\leq 8$.
\end{proposition}
This approach to the problem may even shed light on the validity of Onn's conjecture.\\
The following section outlines the fundamental concepts and provides code examples that were employed to compute the aforementioned data.
\subsection*{Acknowledgements}
    I would like to express my gratitude for the support and guidance provided by my advisor, Uri Onn, during this work as part of the 'research project' course at the Technion in spring 2023. His clear direction and insightful feedback greatly contributed to making this subject accessible. Additionally, I thank Jonathan Shulman for many fruitful discussions.
\begin{center}
\section*{Outline of ideas}
\end{center}
MAGMA's online calculator, available  \href{http://magma.maths.usyd.edu.au/calc/}{here}, fully supports groups of the form $\text{SL}_n(\sfrac{\Z}{p^k})$, and can easily compute their number of conjugacy classes. As an example, the following code computes the number of conjugacy classes for $\text{SL}_4(\sfrac{\Z}{2^3})$:
\begin{Verbatim}
Z:=Integers();
I:=ideal<Z|8>;
A:=quo<Z|I>;
NumberOfClasses(SL(4,A));
\end{Verbatim}
and the output is "1896", as required. Now we return to $\text{SL}_4(\sfrac{\F_2[t]}{(t^3)})$. Unfortunately, MAGMA doesn't support this group. Nevertheless, it does support $\text{GL}_4(\sfrac{\F_2[t]}{(t^3)})$. To compute the number of conjugacy classes of $\text{SL}_4(\sfrac{\F_2[t]}{(t^3)})$ we must find explicit (matrix) generators for $\text{SL}_4(\sfrac{\F_2[t]}{(t^3)})$, as MAGMA can compute the number of conjugacy classes in a group when provided with an explicit set of generators. \\
We begin our search in \cite{PresentationSL3}, where it is shown that the following matrices generate $\text{SL}_3(\sfrac{\Z}{2})$:
\begin{equation*}
    X:=\begin{pmatrix}
        0&1&0\\
        0&0&1\\
        1&0&0
    \end{pmatrix}\;\;\;\; Y:=\begin{pmatrix}
        1&0&1\\
        0&-1&-1\\
        0&1&0
    \end{pmatrix}\;\;\;\;Z:=\begin{pmatrix}
        0&1&0\\
        1&0&0\\
        -1&-1&-1
    \end{pmatrix}
\end{equation*}
So, it seems plausible that the following matrices could potentially generate $SL_4(\sfrac{\Z}{2})$:
\begin{equation*}
    x:=\begin{pmatrix}
        0&1&0&0\\
        0&0&1&0\\
        0&0&0&1\\
        1&0&0&0
    \end{pmatrix}\;\;\;\; y:=\begin{pmatrix}
        1&0&1&0\\
        0&-1&-1&0\\
        0&1&0&0\\
        0&0&0&1
    \end{pmatrix}\;\;\;\;z:=\begin{pmatrix}
        0&1&0&0\\
        1&0&0&0\\
        -1&-1&-1&0\\
        0&0&0&1
    \end{pmatrix}
\end{equation*}
To address the presence of the element $t$, we once again make an informed guess. We include the simplest matrix in $\text{SL}_4(\sfrac{\F_2[t]}{(t^3)})$ that contains $t$:
\begin{equation*}
    w:=\begin{pmatrix}
        1&t&0&0\\
        0&1&0&0\\
        0&0&1&0\\
        0&0&0&1
    \end{pmatrix}
\end{equation*}
Now we compute the size of the subgroup $H$ generated by $\{x,y,z,w\}$ inside $\text{GL}_4(\sfrac{\F_2[t]}{(t^3)})$ and anticipate it is 
the full size of $\text{SL}_4(\sfrac{\F_2[t]}{(t^3)})$, which is\footnote{$|SL_n(\sfrac{\F_p[t]}{(t^k)})|=\frac{p^{(n^2-1)(k-1)}}{p-1}\prod_{0\leq i\leq n-1}(p^n-p^i)$} $21646635171840$. \textbf{If} this is the case, we have successfully found generators for $\text{SL}_4(\sfrac{\F_2[t]}{(t^3)})$ and can thus compute it's number of conjugacy classes. We execute the following code to calculate  $\#H$ and $\#\text{Conj.Classes}(H)$:
\begin{Verbatim}
R<x>:=PolynomialRing(GF(2));
I:=ideal<R|x^3>;
Q:=R/I;
G:=GL(4,Q);
x:=elt<G|
0,1,0,0,
0,0,1,0,
0,0,0,1,
1,0,0,0>;
y:=elt<G|
1,0,1,0,
0,-1,-1,0,
0,1,0,0,
0,0,0,1>;
z:=elt<G|
0,1,0,0,
1,0,0,0,
-1,-1,-1,0,
0,0,0,1>;
w:=elt<G|
1,Q.1,0,0,
0,1,0,0,
0,0,1,0,
0,0,0,1>;
H:=sub<G|x,y,z,w>;
#H;
#ConjugacyClasses(H);
\end{Verbatim}
and the output is $"21646635171840\;\;
1824"$, meaning that $\text{SL}_4(\sfrac{\F_2[t]}{(t^3)})=\langle x,y,z,w\rangle$, and that $\#\text{Conj.Classes}(\text{SL}_4(\sfrac{\F_2[t]}{(t^3)}))=1824$, as required. \\
For the remaining groups, we only provide a set of generators (obtained by similar methods), as the code is very similar.\\
The following matrices generate $\text{SL}_3(\sfrac{\F_3[t]}{(t^3)})$:
\begin{equation*}
    \begin{pmatrix}
        0&1&0\\
        0&0&1\\
        1&0&0
    \end{pmatrix}\;\;\;\; \begin{pmatrix}
        1&0&1\\
        0&-1&-1\\
        0&1&0
 \end{pmatrix}\;\;\;\;\begin{pmatrix}
        0&1&0\\
        1&0&0\\
        -1&-1&-1
    \end{pmatrix}\;\;\;\;\begin{pmatrix}
        1&t&0\\
        0&1&0\\
        0&0&1
    \end{pmatrix}
\end{equation*}
The following matrices generate $\text{SL}_2(\sfrac{\F_2[t]}{(t^k)})$ for $3\leq k\leq 8$:

\begin{equation*}
    \begin{pmatrix}
        1&1\\
        0&1
 \end{pmatrix}\;\;\;\;\begin{pmatrix}
        1&0\\
        t&1
    \end{pmatrix}\;\;\;\;\begin{pmatrix}
        1+t^2&t^2\\
        1&1
    \end{pmatrix}
\end{equation*}
We conclude by noting that we were able to find generators and, consequently, determine the number of conjugacy classes for other groups not mentioned above, including not only $\text{GL}$ but also $\text{SL}$ in various degrees and lengths. For further details and data, interested readers are invited to contact the author.

\bibliographystyle{unsrt}
\bibliography{main}

\end{document}